\documentclass[12pt]{amsart}

\usepackage[latin1]{inputenc}
\usepackage[top=2in, bottom=1.5in, left=1.5in, right=1.5in]{geometry}
\usepackage{amssymb}
\usepackage{graphicx}
\usepackage{amscd}
\usepackage[colorlinks=false,pdfborder=000]{hyperref}
\usepackage{color}
\usepackage{float}
\usepackage{graphics,amsmath,amssymb}
\usepackage{amsthm}
\usepackage{amsfonts}
\usepackage{latexsym}
\usepackage{epsf}
\usepackage{enumerate}

\numberwithin{equation}{section}
\newtheorem{theorem}{Theorem}[section]

\theoremstyle{definition}

\theoremstyle{remark}
\newtheorem{remark}{Remark}[section]

\makeatletter
\def\@settitle{\begin{center}%
  \baselineskip14\p@\relax
  \LARGE
  \@title
  \end{center}%
  \bigskip
}
\makeatother

\begin{document}

\title{A Note on Balancing Binomial Coefficients}

\author{Shane Chern}
\address{School of Mathematical Sciences\\
                Zhejiang University\\
                Hangzhou, 310027, China}
\email{\href{mailto:chenxiaohang92@gmail.com}{chenxiaohang92@gmail.com}}


\textit{Proc. Japan Acad. Ser. A Math. Sci.} \textbf{91} (2015), no. 8, 110--111.

\bigskip \bigskip

\maketitle

\begin{abstract}
In 2014, T. Komatsu and L. Szalay studied the balancing binomial coefficients. In this paper, we focus on the following Diophantine equation
$$\binom{1}{5}+\binom{2}{5}+\cdots+\binom{x-1}{5}=\binom{x+1}{5}+\cdots+\binom{y}{5}$$
where $y>x>5$ are integer unknowns. We prove that the only integral solution is $(x,y)=(14,15)$. Our method is mainly based on the linear form in elliptic logarithms.
\end{abstract}

\section{Introduction}

In a recent paper \cite{Komatsu}, T. Komatsu and L. Szalay studied the balancing binomial coefficients, namely, the Diophantine equation
$$\binom{1}{k}+\binom{2}{k}+\cdots+\binom{x-1}{k}=\binom{x+1}{l}+\cdots+\binom{y}{l}$$
in the positive integer unknowns $x>k$ and $y>x$. In particular, the equation has infinitely many solutions when $k=l=1$ (see also \cite{Finkelstein}), and no solutions when $k=l=2$ or $3$ (see also \cite{Belbachir}). Moreover, when $k=l>3$, it has finitely many solutions. However, for the case $k=l=5$, the authors only got one solution $(x,y)=(14,15)$ through a computer search with $x\le 300$. In this paper, we will completely solve the equation
\begin{equation}\label{eq:1.1}
\binom{1}{5}+\binom{2}{5}+\cdots+\binom{x-1}{5}=\binom{x+1}{5}+\cdots+\binom{y}{5},
\end{equation}
where $y>x>5$, and our result is
\begin{theorem}\label{th:main}
Equation \eqref{eq:1.1} has only one integral solution $(x,y)=(14,15)$.
\end{theorem}
Our method of proof, which is mainly based on the linear form in elliptic logarithms, is motivated by \cite{Ingram}, and further by \cite{Stroeker}. Throughout this paper, we use the notations of \cite{Stroeker}.

\section{Proof of Theorem \ref{th:main}}

Note that
$$\binom{1}{k}+\binom{2}{k}+\cdots+\binom{x}{k}=\binom{x+1}{k+1},$$
we can therefore rewrite \eqref{eq:1.1} as
\begin{equation}\label{eq:2.1}
\binom{x}{6}+\binom{x+1}{6}=\binom{y+1}{6}.
\end{equation}
Set $u=(x-2)^2$ and $v=(y-1)(y-2)$, \eqref{eq:2.1} becomes
\begin{equation}\label{eq:2.2}
2u^3-10u^2+8u=v^3-8v^2+12v.
\end{equation}
The transformation
$$X=\dfrac{-4u-17v+58}{2u-3v},$$
$$Y=\dfrac{-146u^2-5v^2+686u-188v}{(2u-3v)^2}$$
yields a minimal Weierstrass model for \eqref{eq:2.2}, specifically,
\begin{equation}\label{eq:2.3}
E:Y^2=X^3-X^2-30X+81.
\end{equation}
With \textit{Magma}, we can find the Mordell-Weil group $E(\mathbb{Q})\cong\mathbb{Z}\oplus\mathbb{Z}\oplus\mathbb{Z}$ is generated by the points $P_1=(3,-3)$, $P_2=(-6,3)$ and $P_3=(11,31)$.

Let $Q_0=(X_0,Y_0)$ be the image on \eqref{eq:2.3}  of the point at infinity on \eqref{eq:2.2}, then we have $X_0=7+2\alpha+3\alpha^2$ and $Y_0=-17-15\alpha-8\alpha^2$ where $\alpha=\sqrt[3]{2}$. Note that $Q_0\in E(\mathbb{Q}(\alpha))$. Also note that $v=\alpha u+\beta$ where $\beta=(8-5\alpha)/3$ is the asymptote of the curve.

It is easy to verify that
\begin{equation}\label{eq:2.4}
\dfrac{dv}{6u^2-20u+8}=-\dfrac{1}{4}\dfrac{dX}{Y}.
\end{equation}
Note that for $v\ge 30$, $u(v)$ given by \eqref{eq:2.2} can be viewed as a strictly increasing function of $v$, we therefore have
\begin{equation}\label{eq:2.5}
\int_v^\infty\dfrac{dv}{6u^2-20u+8}=\dfrac{1}{4}\int_{X_0}^X\dfrac{dX}{Y}.
\end{equation}
It is also easy to verify that $6u^2-20u+8>3v^2$ for $v\ge 30$. Thus, we have
\begin{equation}\label{eq:2.6}
\int_v^\infty\dfrac{dv}{6u^2-20u+8}<\dfrac{1}{3}\int_v^\infty\dfrac{dv}{v^2}=\dfrac{1}{3v}.
\end{equation}

Let $P=m_1P_1+m_2P_2+m_3P_3$ be an arbitrary point on \eqref{eq:2.3} with integral coordinates $u$, $v$ on \eqref{eq:2.2}. We have
\begin{equation}\label{eq:2.7}
\int_{X_0}^X\dfrac{dX}{Y}=\int_{X_0}^\infty\dfrac{dX}{Y}-\int_{X}^\infty\dfrac{dX}{Y}=\omega(\phi(Q_0)-\phi(P)),
\end{equation}
where $\omega=5.832948\ldots$ is the fundamental real period of $E$, and
\begin{align*}
\phi(P)&=\phi(m_1P_1+m_2P_2+m_3P_3)\\
&=m_1\phi(P_1)+m_2\phi(P_2)+m_3\phi(P_3)+m_0
\end{align*}
with $m_0\in\mathbb{Z}$ and all $\phi$-function are in $[0,1)$. Put $M=\max_{1\le i\le3}|m_i|$, it follows $|m_0|\le 3M$. By Zagier's algorithm (see \cite{Zagier}), we have $u_1=\omega\phi(P_1)=4.158074\ldots$, $u_2=\omega\phi(P_2)=2.851605\ldots$, $u_3=\omega\phi(P_3)=0.627538\ldots$, and $u_0=\omega\phi(Q_0)=5.289657\ldots$. Let
\begin{align*}
L(P)&=\omega(\phi(Q_0)-\phi(P))\\
&=u_0-m_0\omega-m_1u_1-m_2u_2-m_3u_3,
\end{align*}
we then obtain the lower bound
\begin{equation}\label{eq:lb}
|L(P)|>\exp(-c_4(\log(3M)+c_5)(\log\log(3M)+c_6)^6),
\end{equation}
where $c_4=7\times10^{160}$, $c_5=2.1$, and $c_6=21.2$, by applying David's result \cite{David} (see also \cite{Tzanakis}).

By \eqref{eq:2.5} and \eqref{eq:2.6}, we also have
\begin{equation}\label{eq:2.8}
|L(P)|=4\int_v^\infty\dfrac{dv}{6u^2-20u+8}\le\dfrac{4}{3v}.
\end{equation}
For $v\ge30$, it is easy to verify that
\begin{equation}\label{eq:2.9}
h(P)\le\log(4u+17v-58)<3.044523+\log v.
\end{equation}
Here $u$ and $v$ are required to be integral. Moreover we have
\begin{equation}\label{eq:2.10}
\hat{h}(P)\ge c_1M^2
\end{equation}
where $c_1=0.125612\ldots$ is the least eigenvalue of the N\'eron-Tate height pairing matrix. Note that Silverman's bound for the difference of heights on elliptic curves gives that
\begin{equation}\label{eq:2.11}
2\hat{h}(P)-h(P)<7.846685.
\end{equation}
By \eqref{eq:2.8}, \eqref{eq:2.9}, \eqref{eq:2.10}, and \eqref{eq:2.11}, we obtain the upper bound
\begin{equation}\label{eq:ub}
|L(P)|<\exp(11.1789-0.251224M^2).
\end{equation}

Together with \eqref{eq:lb} and \eqref{eq:ub}, we therefore find an absolute upper bound $M_0=1.4\times10^{86}$ for $M$. Applying the LLL algorithm (cf. \cite{Stroeker}), we may reduce this bound to $M=11$. Through a computer search, we therefore prove that $(x,y)=(14,15)$ is the only integral solution of \eqref{eq:1.1}.

\begin{remark}
It is of interest to mention that all integral solutions $(u,v)$ of \eqref{eq:2.2} are given in Table \ref{ta:01}, by slightly modifying our proof and then through a similar computer search.
\end{remark}

\begin{table}[ht]
\caption{Integral solutions $(u,v)$ of \eqref{eq:2.2}}\label{ta:01}
\renewcommand\arraystretch{1.25}
\noindent\[
\begin{tabular}{p{0.1cm} r@{,}l p{0.2cm} r@{,}l p{0.2cm} r@{,}l p{0.1cm}}
\hline
& $(0$&$0)$ & & $(0$&$2)$ & & $(0$&$6)$ &\\
& $(1$&$0)$ & & $(1$&$2)$ & & $(1$&$6)$ &\\
& $(4$&$0)$ & & $(4$&$2)$ & & $(4$&$6)$ &\\
& $(9$&$12)$ & & $(144$&$182)$ & & $(-56$&$-70)$ &\\
\hline
\end{tabular}
\]
\end{table}

\bibliographystyle{amsplain}

\end{document}